\documentclass[10pt,a4paper,preprint,review]{elsarticle}
\usepackage[top=1 in, bottom=1 in, left=1 in, right=1 in, letterpaper]{geometry}

\usepackage[english]{babel}

\usepackage{amsmath}
\usepackage{amssymb}
\usepackage{calc,comment}

\usepackage{tikz}

\usetikzlibrary{arrows,arrows.meta}
\usetikzlibrary{calc}

\pgfdeclarelayer{bg}
\pgfsetlayers{bg,main}

\newtheorem{example}{Example}

\newtheorem{definition}{Definition}

\usepackage{mathtools}
\newcommand\pnpequality{\stackrel{\mathclap{\normalfont\mbox{\footnotesize?}}}{=}}

\usepackage[capitalize]{cleveref}

\usepackage{multicol, setspace}
\usepackage{multirow}
\usepackage{dashrule}

\usepackage{booktabs}

\usepackage{longtable, tabu}
\usepackage{ltablex}

\usepackage{caption}
\usepackage{subcaption}

\newcommand{\xil}{x_{i, \, l}}
\newcommand{\ylj}{y_{l, \, j}}
\newcommand{\zij}{z_{i, \, j}}

\newcommand{\plk}{p_{l, \, k}}
\newcommand{\qil}{q_{i, \, l}}
\newcommand{\ci}{c_i}

\newcommand{\Ail}{A_i^L}
\newcommand{\Aiu}{A_i^U}
\newcommand{\Sl}{S_l}
\newcommand{\Djl}{D_j^L}
\newcommand{\Dju}{D_j^U}
\newcommand{\Cik}{C_{i, \, k}}

\newcommand{\Pjkl}{P_{j, \, k}^L}
\newcommand{\Pjku}{P_{j, \, k}^U}

\newcommand{\vilj}{v_{i, \, l, \, j}}

\newcommand{\vect}[1]{\boldsymbol{#1}}

\newcommand{\pij}{p_{i,j}}
\newcommand{\xijt}{x_{i,j,t}}
\newcommand{\bijl}{b_{i,j}^L}
\newcommand{\biju}{b_{i,j}^U}
\newcommand{\bijt}{b_{i,j,t}}
\newcommand{\yst}{y_{s,t}}

\newcommand{\qism}{q_{i,s}^-}

\newcommand{\stin}{\text{in}}
\newcommand{\stout}{\text{out}}

\journal{}

\begin{document}

\begin{frontmatter}


\title{ \textbf{Approximation Algorithms for Process Systems Engineering} }

\author{Dimitrios Letsios}
\author{Radu Baltean-Lugojan}
\author{Francesco Ceccon}
\author{Miten Mistry}
\author{Johannes Wiebe}
\author{Ruth Misener}

\address{Department of Computing, Imperial College London, SW7 2AZ, UK}

\ead{r.misener@imperial.ac.uk}

\date{}

\doublespacing

\begin{abstract}
Designing and analyzing algorithms with provable performance guarantees enables 
efficient optimization problem solving in different application domains, e.g.\ communication networks, transportation, economics, and manufacturing.
Despite the significant contributions of approximation algorithms in engineering, only limited and isolated works contribute from this perspective in process systems engineering.
The current paper discusses three representative, $\mathcal{NP}$-hard problems in process systems engineering: (i) pooling, (ii) process scheduling, and (iii) heat exchanger network synthesis. 
We survey relevant results and raise major open questions.
Further, we present approximation algorithms applications which are relevant to process systems engineering: (i) better mathematical modeling, (ii) problem classification, (iii) designing solution methods, and (iv) dealing with uncertainty.
This paper aims to motivate further research at the intersection of approximation algorithms and process systems engineering.
\end{abstract}

\end{frontmatter}



\section{Introduction}
Both Theoretical Computer Science (TCS) and Process Systems Engineering (PSE) design efficient algorithms for challenging optimization problems. But, although both areas consider mixed-integer [non-]linear optimization, the two domains have several intellectual differences:

\begin{enumerate}
\item TCS often addresses worst-case instances whereas PSE typically solves challenging, practically-relevant instances,
\item TCS analytically derives theoretical performance guarantees while PSE computationally proves global optimality, e.g.\ using a mixed-integer [non-]linear optimization solver,
\item TCS uses bottom-up approaches, e.g.\ solving interesting problem special cases and determining polynomiality and inapproximability boundaries. Meanwhile, PSE more often employs top-down techniques, e.g.\ problem decompositions such as Dantzig-Wolfe or Benders, 
\item TCS frequently concerns itself with near-optimal algorithms giving a provably-good feasible solution while PSE is more interested in the deterministic global solution,
\item TCS focusses on polynomial tractability whereas PSE is more interested in practical computational scalability for industrial instances.
\end{enumerate}


The design and analysis of \emph{approximation algorithms}, i.e.\ heuristics with performance guarantees, is well-established in TCS \citep{Hochbaum1996,Schulz1997,Vazirani2001,Williamson2011}.
\citet{Johnson1974} introduces approximation algorithms as a general framework for solving combinatorial optimization problems.
\citet{Sahni1977} presents an early tutorial of techniques.
Approximation algorithms compute feasible solutions that are provably close to optimal solutions. These approximation algorithms are designed to have efficient, i.e.\ polynomial, running times. Approximation algorithms have been developed for optimization problems arising in application domains,
including 
communication networks \citep{Garg1996,Goemans1994,Johnson1978,Kleinberg1996,Leighton1999},
transportation  \citep{Christofides1976,Frederickson1976,Golden1980,Kruskal1956,Laporte1992,Rosenkrantz1977}, 
economics \citep{Daskalakis2006,Daskalakis2009,Lipton2003,Papadimitriou2014}, 
and manufacturing \citep{Gonzalez1978,Graham1969,Hall1989,Jackson1955}.
But approximation algorithms have not received much attention in PSE.
The PSE community is mainly interested in global optimization methods because suboptimal solutions may incur significant costs, or even be incorrect \citep{Grossmann2013}.
At a first glance, approximation algorithms do not fit the PSE preference towards an exact solution.
Furthermore, heuristics with performance guarantees cannot fully address the very complex, highly inapproximable, industrially-relevant optimization problems in PSE.

\begin{quote}
This paper argues that, contrary to the aforementioned, surface-level distinctions, approximation algorithms are deeply applicable to PSE. 
We substantiate our claims by offering applications where approximation algorithms can be particularly useful for solving challenging process systems engineering optimization problems.
\end{quote}

In the last 30 years, there has been significant progress in designing approximation algorithms and understanding the limits of proving analytical performance guarantees.
Problems that sound simple, e.g.\ makespan scheduling and bin packing, are believed to be hard \citep{Garey2002}.
Computational complexity theory and $\mathcal{NP}$-hardness provide a mathematical foundation for this belief.
Under the widely adopted conjecture $\mathcal{P}\neq\mathcal{NP}$, no algorithm can solve an $\mathcal{NP}$-hard problem in polynomial worst-case running time, e.g.\ scheduling $n$ jobs in time proportional to some polynomial $p(n)$ of $n$.
Approximation algorithms cope with $\mathcal{NP}$-hardness by producing, in polynomial time, good suboptimal solutions. 
In particular, a $\rho$-approximate algorithm for a minimization (resp.\ maximization) problem computes, for every input, a solution of cost (resp.\ profit) at most (resp.\ least) $\rho$ times the optimum.
The performance guarantee $\rho$ quantifies the worst-case distance of an approximation algorithm's solutions from being optimal, i.e.\
provides an optimality gap for 
pathological optimization problem instances. 
But, in practice, an approximation algorithm may produce a significantly better solution than the worst-case bound.
From a complementary viewpoint, $\mathcal{NP}$-hardness specifies the limits in developing optimization algorithms with polynomial worst-case running times. 
Meanwhile, hardness of approximation settles the limits of polynomial approximation algorithms, e.g.\ may prove that designing $\rho$-approximation algorithm with small $\rho$ is impossible. 


The manuscript proceeds as follows:
Section~\ref{Section:Approximation_Algorithms} introduces approximation algorithms.
Section~\ref{Section:Applications} presents applications of approximation algorithms in PSE.
The remainder of the paper provides a brief survey of three major PSE optimization problems: Sections~\ref{Section:Pooling}-
\ref{Section:Heat_Exchangers} discuss pooling, process scheduling, and heat exchanger network synthesis, respectively. 
These sections present a collection of $\mathcal{NP}$-hard problems for which approximation algorithms can be useful.
Section~\ref{Section:Conclusion} concludes the paper.

\section{Approximation Algorithms}
\label{Section:Approximation_Algorithms}

This section introduces the notion of an \emph{approximation algorithm}, i.e.\ a heuristic with a performance guarantee \citep{Vazirani2001, Williamson2011}.
Many optimization problems are $\mathcal{NP}$-hard and,
under the widely adopted conjecture $\mathcal{P}\neq\mathcal{NP}$, 
there is no polynomial algorithm solving an $\mathcal{NP}$-hard problem.
Approximation algorithms investigate the trade-off between optimality and computational efficiency for a range of applications. 

An approximation algorithm is a polynomial algorithm producing a near-optimal solution to an optimization problem. 
Formally, consider an optimization problem, without loss of generality a minimization problem, and a polynomial Algorithm $A$ for getting a feasible solution (not necessarily the global optimum).

\begin{definition}[\citep{Johnson1974}]
\label{Def:Approximation_Algorithm}
For each problem instance $I$, let $C_{A}(I)$ and $C_{OPT}(I)$ be the algorithm's objective value and the globally optimal objective value, respectively.
Algorithm $A$ is $\rho$-approximate if, for every problem instance $I$, it holds that
$C_{A}(I)\leq \rho\cdot C_{OPT}(I)$.
Value $\rho$ is the \emph{approximation ratio} of Algorithm $A$. 
\end{definition}

That is, a $\rho$-approximation algorithm computes, in polynomial time, a solution with an objective value at most $\rho$ times the optimal objective value. 
Since $C_{OPT}(I)$ is the globally optimal objective value, we trivially note that $\rho > 1$. Theoretical computer scientists may seek constant $\rho$, e.g.\ $\rho = \frac{4}{3}$, since (i) the algorithm will not degrade with growing problem size and (ii) a 
constant approximation ratio is significant progress in solving an optimization problem efficiently. 

\begin{figure}
\begin{center}
\includegraphics[]{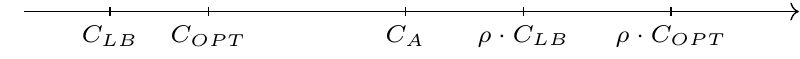}
\caption{Analysis of an Approximation Algorithm}
\label{Fig:ApproximationAlgorithm}
\end{center}
\end{figure}

In general, to prove a $\rho$-approximation ratio, we proceed as depicted in Figure \ref{Fig:ApproximationAlgorithm}.
For each problem instance, we compute analytically a lower bound $C_{LB}(I)$ of the optimal objective value, i.e.\ $C_{LB}(I) \leq C_{OPT}(I)$. One lower bounding method, replacing $\{0, 1\}$ binary variables with a fractional $[0,1]$ relaxation, will be familiar to the PSE community from mixed-integer optimization. 
Other common lower bounding methods include basic packing and covering \citep{Chvatal1979, Graham1969, Mcnaughton1959}, duality \citep{Cornuejols1977, Held1970}, and semidefinite programming \citep{Goemans1997,Goemans1995}.
%
The next step proves that an algorithm's objective value is at most $\rho$ times the lower bound, i.e.\ $C_{A}(I)\leq \rho \cdot C_{LB}(I)$.
Therefore, proving a $\rho$-approximation is, in some sense, equivalent to matching an upper objective bound with a lower objective bound.
A $\rho$-approximation ratio is \emph{tight} for Algorithm $A$, if we can prove that there is no lower ratio.


Definition~\ref{Def:Approximation_Scheme} introduces a well-known family of algorithms known as \emph{approximation schemes}. 

\begin{definition}[\citep{Sahni1976}]
\label{Def:Approximation_Scheme}
Algorithm $A$ is an approximation scheme if, for every problem instance $I$ and input parameter $\epsilon>0$, it holds that $C_A(I)\leq(1+\epsilon)C_{OPT}(I)$\footnote{For maximization problems, the performance bound becomes $C_A(I)\geq(1-\epsilon)C_{OPT}(I)$.}.
If the running time of $A$ is bounded by a polynomial in the instance size, then $A$ is a \emph{polynomial-time approximation scheme (PTAS)}.
When $A$ is also polynomial in $1/\epsilon$, then it is called a \emph{fully polynomial-time approximation scheme (FPTAS)}.
\end{definition}

A PTAS, and in particular a FPTAS, is the best approximation result that one can hope for an $\mathcal{NP}$-hard problem, unless $\mathcal{P}=\mathcal{NP}$.
\citet{Schuurman2000} provide a tutorial 
for designing PTAS.
Despite their theoretical significance, PTAS are often not necessarily the most competitive algorithms for real-world problem instances, e.g.\ in the Euclidean traveling salesman problem case \citep{Johnson2012,Johnson1997}.

Obtaining PTAS or even constant performance guarantee in polynomial time might not be possible. 
Hardness of approximation provides a toolbox of techniques for deriving such negative results \citep{Goderbauer2019}.
The best possible performance ratios for problems that do not admit a PTAS are typically $O(1)$, $O(\log n)$, or $O(n^{O(1)})$.
The standard \emph{big-O notation} $O(\cdot)$ means that $O(1)$ is some constant while $O(\log n)$ and $O(n^{O(1)})$ indicate some function of $n$ asymptotically upper bounded by a logarithmic and polynomial function, respectively.


\section{Applications of Approximation Algorithms in PSE}
\label{Section:Applications}

This section discusses ways of using approximation algorithms in PSE and presents examples of past contributions motivating research in these directions. 




\subsection{Mathematical Modeling}

PSE optimizes chemical, biological, and physical processes using systematic computer-aided approaches.
A significant part of the PSE literature is devoted to evaluating, verifying, refining, and validating mathematical models capturing natural phenomena. 
These models quantitatively predict process outputs subject to initial conditions.
Frequently, dealing with PSE problems involves solving mixed-integer linear programming (MILP) models.
Modern MILP solver performance depends considerably on the underlying MILP formulation \citep{Vielma2015}.
Critical MILP formulation aspects include the size and strength of the LP relaxation.
The theory of approximation algorithms provides a methodology for evaluating relaxation quality using worst-case analysis.
Designing strong relaxations with analytically proven performance guarantees (i) reveals meaningful insights for relaxations that are well-suited for a MILP instance and (ii) derives effective reformulations towards those structures \citep{Cornuejols2008}.
Structural combinatorial properties of near-optimal solutions may strengthen MILP models, e.g.\ with valid inequalities and symmetry breaking constraints \citep{Margot2010}.

\begin{example}
Tight theoretical bounds show that the so-called $PQ$-formulation attains the best possible performance guarantee \citep{Dey2015,Tawarmalani2002}.
Empirical evidence demonstrates the computational superiority of the $PQ$-formulation compared to other formulations 
\citep{Alfaki2013multicommodity}.
\end{example}


\subsection{Problem Classification}

PSE investigates different techniques, e.g.\ branch-and-bound, cutting planes, metaheuristics, and decompositions, for effectively solving a variety of MILP problems arising in engineering.
A major goal is developing general purpose solvers 
selecting the most appropriate solution method for each concrete MILP instance.
Heuristics are an essential tool in MILP solving \citep{Berthold2014, Fischetti2010, Schulz1997, Williamson2011}.
Performance guarantees, which arise from the theory of approximation algorithms, 
evaluate and classify the computational performance of heuristics.
Investigating the approximation properties of $\mathcal{NP}$-hard problems involves (i) designing efficient approximation algorithms and (ii) determining inapproximability results.
Approximation algorithms exploit special structure and expose tractable optimization problem subcases.
Hardness of approximation classifies problems
from a computational complexity viewpoint and determines the limits of efficient approximation.
These directions contribute to selecting and employing suitable solution methods for PSE problems.


\begin{example}
\citet{Chen1998} and \citet{Coffman2013} provide extensive reviews and classifications of approximation algorithms for scheduling and bin packing problems.
These algorithm portfolios improve the ability of solving such problems efficiently \citep{Bischl2016}.
\end{example}





\subsection{Design of Solution Methods}

PSE optimization methods can be broadly divided into global and local (nonlinear programming) \citep{Grossmann2013}.
Global optimization has attracted substantial attention by the PSE community because it overcomes the limitations of local optimization in generating solutions with guarantees of $\epsilon$-global optimality.
On the other hand, local optimization can be particularly useful when dealing with PSE problems of massive size.
TCS provides a framework for designing algorithms attaining good trade-offs in terms of solution quality and running time efficiency \citep{Schulz1997}.
An approximation algorithm computes solutions quickly that are provably close to optimal.
Furthermore, TCS offers a toolbox of techniques for designing approximation algorithms including local search, dynamic programming, linear programming, duality, semidefinite programming, and randomization.
Hence, approximation algorithms are handy for very large-scale PSE problem solving with certified distance from optimality.

\begin{example}
\citet{Letsios2018} provide a collection of heuristics with proven performance guarantees for solving large-scale instances of the minimum number of matches problem in heat recovery network design.
These heuristics obtain better solutions than commercial solvers in reasonable time frames.
\end{example}

\subsection{Dealing with Uncertainty}

Process operations exhibit inherent uncertainty such as demand fluctuations, equipment failures, and temperature variations.
The successful application of PSE optimization models in practice depends crucially on the ability to handle uncertainty \citep{PISTIKOPOULOS1995553}.
A key challenge is to construct robust solutions and determine suitable recovery actions responding proactively and reactively to variations and unexpected events. 
Optimization methods under uncertainty yield suboptimal solutions with respect to the ones that may be obtained with full input knowledge.
Approximate performance guarantees are useful for characterizing the structure of robust solutions \citep{Bertsimas2004, Bertsimas2011, Goerigk2016}.
Recovery strategies improve robust solutions by making second-stage decisions after the uncertainty is revealed \citep{Liebchen2009}.
TCS approaches derive recovery methods attaining good trade-offs in terms of final solution quality and initial solution transformation cost \citep{Ausiello2011, Schieber2018, Skutella2016}.


\begin{example}
Past literature obtains useful structural properties of robust solutions for fundamental combinatorial optimization problems under uncertainty.
\citet{Monaci2013} show that the number of perturbed item weights does not affect the solution quality in the knapsack problem.
\citet{Letsios2018scheduling} show that lexicographic optimization imposes optimal substructure for the makespan scheduling problem.
\citet{Schieber2018} present a framework designing reoptimization algorithms with analytically proven performance guarantees and present a family of fully polynomial-time reapproximation schemes.
\end{example}

\section{Pooling}
\label{Section:Pooling}

Pooling is a major optimization problem with applications, e.g.\ in petroleum refining \citep{Baker1985}, crude oil scheduling \citep{Lee1996, Li2012}, natural gas production \citep{Li2011,Selot2008}, hybrid energy systems \citep{Baliban2012}, water networks \citep{Galan1998}, and a sub-problem in general mixed-integer nonlinear programs (MINLP) \citep{Ceccon2016}.
The goal is to blend raw materials in intermediate pools in order to produce final products, minimizing process costs while satisfying customer demand and meeting final product requirements.
Pooling is an $\mathcal{NP}$-hard, nonconvex nonlinear optimization problem (NLP) and variant of network flow problems.
The challenge is to deal with bilinear terms and multiple local minima \citep{Haverly1978}.



\subsection{Brief Literature Overview}

Algorithms for the pooling problem have evolved in tandem with state-of-the-art non-convex quadratically-constrained optimization solvers \citep{Audet2004, Boukouvala2016, Misener2009}. 
Early approaches rely on sparsity and tackle large-scale instances with successive linear programming (SLP), i.e.\ efficiently solving a sequence of linear programs obtained by first-order Taylor approximations of bilinear terms \citep{Baker1985, Dewitt1989}. 
\cite{Visweswaran1990, Visweswaran1993} 
investigate global optimization methods using duality theory and Lagrangian relaxations, which are further explored by \citet{Adhya1999}. 
A subsequent line of work develops strong relaxations with convex envelopes including 
reformulation-linearization cuts \citep{Meyer2006, Quesada1995, Sherali1999, Sherali1992}, 
McCormick envelopes \citep{AlKhayyal1983, Foulds1992, McCormick1976}, 
sum-of-squares \citep{Marandi2018},
multi-term and edge concave cuts \citep{Bao2009, Misener2011miqcqp, Misener2012cuts}. 
These approaches are employed in state-of-the-art mixed-integer nonlinear programming software where piecewise-linear relaxations may further improve solver performance on pooling problems \citep{Gounaris2009, Hasan2010, Kolodziej2013, Misener2011pooling, Misener2011miqcqp, Wicaksono2008}. Further valid linear and convex inequalities are derived from nonconvex restrictions of the pooling problem \citep{luedtke2018strong}. Parametric uncertainty in the pooling problem has recently been considered using stochastic programming and robust optimization approaches \citep{Li2012,Li2011,Wiebe2019}.

Exact MINLP methods exhibit exponential worst-case behavior,
so designing heuristic approaches with analytically proven performance guarantees is useful for (i) finding provably good solutions on a fast time frame and (ii) solving very large scale instances where.

\subsection{Problem Definition}
\label{Section:Pooling_Problem_Definition}

\begin{figure}
    \centering
    \includegraphics[width=0.7\textwidth]{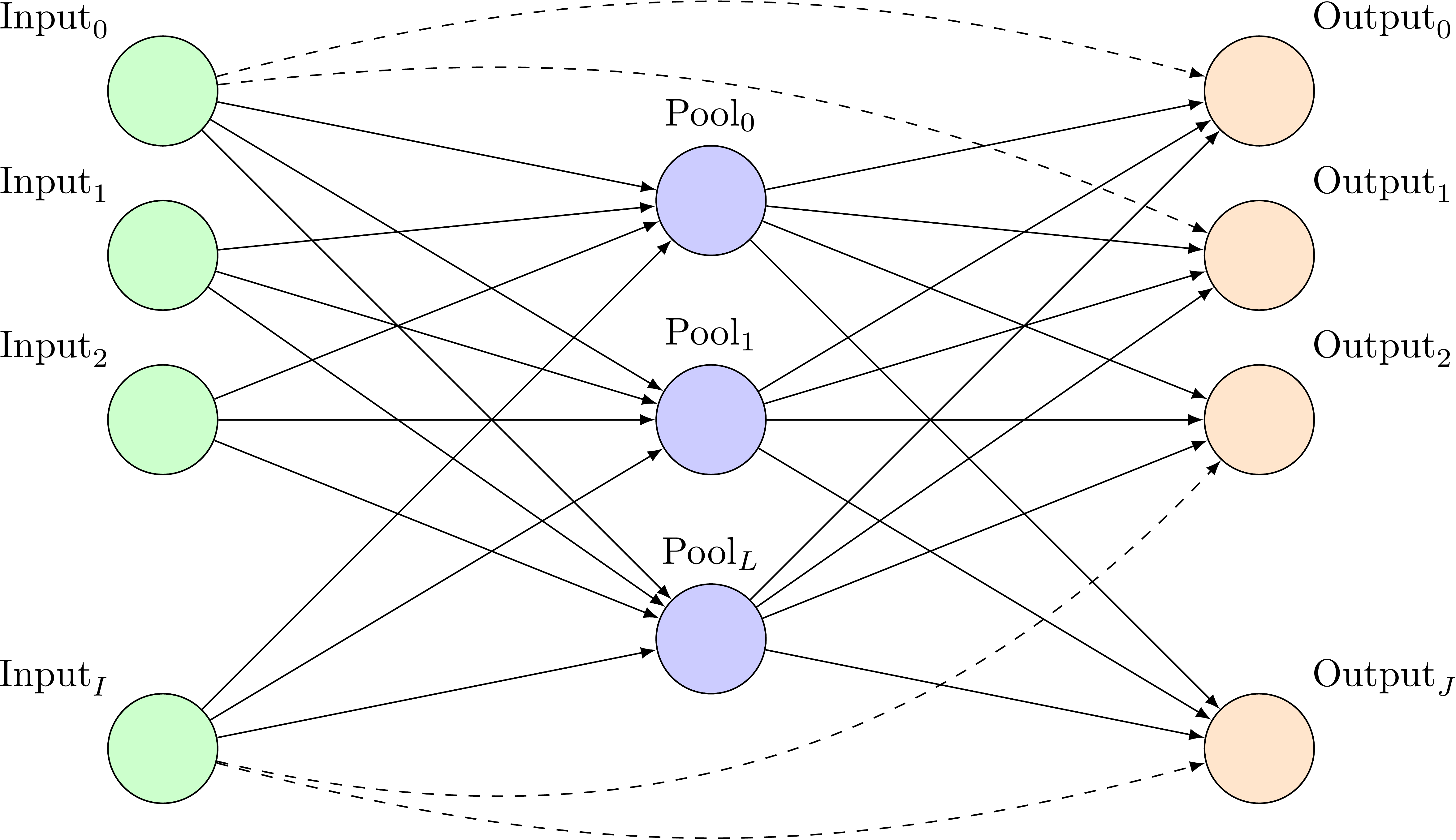}
    \caption{Pooling network $T=(N,A)$ with (i) input, (ii) pool, and (iii) output nodes. Straight arcs represent input-to-pool and pool-to-output flow. Dashed arcs illustrate bypass input-to-output flow.}
    \label{fig:pooling-network}
\end{figure}

A pooling problem instance is a directed network $T=(N,A)$, where $N$ is the set of vertices and $A$ is the set of arcs.
Figure~\ref{fig:pooling-network} illustrates a pooling network.
The set of nodes can be partitioned into the sets $I\cup L\cup J$, where $I$ is the set of \emph{input} or \emph{source nodes}, $L$ is the set of \emph{pool nodes}, and $J$ is the set of \emph{output} or \emph{terminal nodes}. 
The directed arcs $A$ are a subset of $X\cup Y\cup Z$, where $X=I\times L$, $Y=L\times J$, and $Z=I\times J$.
A solution to the pooling problem can be viewed as a flow of materials in the network $T$.
Input nodes introduce raw materials, pool nodes mix raw materials and produce intermediate products, while output nodes export final products.
Let $x_{i,l}$ be the flow exiting input node $i\in I$ and entering pool $l\in L$.
Then, $\sum_{i\in I}x_{i,l}$ units of flow enter pool $l\in L$.
Similarly, denote by $y_{l,j}$ and $z_{i,j}$ the flow transferred from pool $l\in L$ to output node $j\in J$ and the bypass flow routed directly from input node $i\in I$ to terminal node $j\in J$, respectively.
Then, $\sum_{l\in L}y_{l,j} + \sum_{i\in I}z_{i,j}$ units of flow enter output node $j\in J$.
Flow conservation enforces that the amount $\sum_{i\in I}x_{i,l}$ of entering raw materials is equal to the amount $\sum_{j\in J}y_{l,j}$ of exiting intermediate product, for each pool $l\in L$.
The total quantity of raw material $i\in I$ and final product $j\in J$ are subject to the box constraints
$\Ail\leq \sum_{i\in I}x_{i,l}\leq \Aiu$ and $\Djl\leq \sum_{l\in L}y_{l,j} + \sum_{i\in I}z_{i,j}\leq\Dju$, respectively.
Moreover, pool $l\in L$ has flow capacity $S_l$. 
\ref{app:pooling_notation} presents the notation for the pooling problem.

Pooling problem monitors a set $K$ of quality attributes, e.g.\ concentrations of different chemicals, for each raw material, intermediate, and final product.
Raw material $i\in I$ has attribute $k\in K$ value $\Cik$.
The intermediate and final product attribute values are determined assuming linear blending.
Specifically, the attribute $k\in K$ value $p_{l,k}$ in pool $l\in L$ satisfies $p_{l,k} \sum_{j\in J}y_{l,j} = \sum_{i\in I} x_{i,l}C_{i,k}$.
On the other hand, the attribute $k\in K$ value of end product $j\in J$ is $\sum_{l\in L}p_{l,k}y_{l,j} + \sum_{i\in I}z_{i,j}C_{i,k}$.
Final product $j\in J$ is constrained to admit attribute $k\in K$ value in the range $[\Pjkl, \Pjku]$.
The goal is to optimize raw material costs and sales profit.
In particular, let $c_i$ and $d_j$ be the unitary cost of raw material $i$ and the unitary profit of end product $j$.
Then, the pooling problem minimizes $\sum_{i\in I}c_ix_i - \sum_{j\in J}d_jy_j$.



\subsection{Mathematical Models}


This section provides the standard $P$- and $PQ$-formulations \citep{Ben-Tal1994,Haverly1978,Quesada1995,Tawarmalani2002} for modeling the pooling problem.
For simplicity, these formulations are presented assuming the pooling network is complete, i.e.\ contains all possible arcs, but can be easily extended to arbitrary networks.

\subsubsection{$P$-formulation}

The $P$-formulation \citep{Haverly1978} uses the Section~\ref{Section:Pooling_Problem_Definition} arc flow and intermediate product attribute variables and can be stated using Equations~(\ref{Eq:P_Model}). 
The resulting quadratic programming formulation includes bilinear terms due to linear blending. 

{ \footnotesize
\begin{subequations}
\label{Eq:P_Model}
\begin{align}
\min_{\vect{x}, \vect{y}, \vect{z}, \vect{p}} &
\sum_{(i, l) \in X} \ci \xil - \sum_{(l, j) \in Y} d_j \ylj - \sum_{(i, j) \in Z} (d_j - \ci) \zij \label{Eq:P_Model_Objective} \\
& \Ail \leq \sum\limits_{l \in L} \xil + \sum_{j \in J} \zij \leq  \Aiu & i\in I \label{Eq:P_Model_Supply_Bounds} \\
& \sum\limits_{j \in J} \ylj \leq \Sl & l\in L \label{Eq:P_Model_Pool_Capacity} \\
& \Djl \leq \sum\limits_{l \in L} \ylj + \sum_{i \in I} \zij \leq \Dju & j\in J \label{Eq:P_Model_Demand_Bounds} \\
& \sum_{i \in I} \xil = \sum_{j \in J} \ylj & l\in L \label{Eq:P_Model_Pool_Conservation} \\
& \sum_{i \in I} \Cik \xil = \plk \sum_{j \in J} \ylj & l\in L, k\in K \label{Eq:P_Model_Quality_Balances} \\
& \sum_{l \in L} \plk \ylj + \sum_{i \in I} \Cik \zij \geq \Pjkl ( \sum_{l \in L} \ylj + \sum_{i \in I} \zij) & j\in J, k\in K 
\label{Eq:P_Model_Quality_Lower} \\
& \sum_{l \in L} \plk \ylj + \sum_{i \in I} \Cik \zij \leq \Pjku ( \sum_{l \in l} \ylj + \sum_{i \in I} \zij) & j\in J, k\in K 
\label{Eq:P_Model_Quality_Upper} \\
& \xil, \ylj, \zij \geq 0 & i\in I, l\in L, j\in J \label{Eq:P_Model_Flow_Bounds} \\
& \plk \geq 0 & l\in L, k\in K \label{Eq:P_Model_Quality_Bound}  
\end{align}
\end{subequations}
}
Expression (\ref{Eq:P_Model_Objective}) optimizes raw material costs and final product profits.
Constraints (\ref{Eq:P_Model_Supply_Bounds}) - (\ref{Eq:P_Model_Demand_Bounds}) impose bounds on the raw material quantities, pool sizes, and final product quantities, respectively.
Constraints (\ref{Eq:P_Model_Pool_Conservation}) ensure material balances.
Constraints (\ref{Eq:P_Model_Quality_Balances}) model linear blending.
Constraints (\ref{Eq:P_Model_Quality_Lower}) - (\ref{Eq:P_Model_Quality_Upper}) enforce quality specifications.
Constraints (\ref{Eq:P_Model_Flow_Bounds}) - (\ref{Eq:P_Model_Quality_Bound}) ensure that all variables are non-negative.


\subsubsection{$PQ$-formulation}

The $PQ$-formulation \citep{Quesada1995,Tawarmalani2002} extends the \citet{Ben-Tal1994} pooling problem $Q$-formulation and replaces the $P$-formulation variables $\xil$ with path flow variables $\vilj$ and proportion variables $\qil$, for $i\in I$, $l\in L$, and $j\in J$. 
Specifically, $\vilj$ represents the flow transferred from input node $i\in i$ to output node $j\in J$ via pool node $l\in L$ and $\qil$ corresponds to the proportion of total flow entering pool $l\in L$ that originates from input node $i\in I$, i.e.\ $\xil = \qil \sum_{j \in J} \ylj$ and $0\leq \qil\leq 1$.
The $PQ$-formulation can be stated using Equations (\ref{Eq:PQ_Model}) and results in a tighter McCormick relaxation compared to the $P$-formulation. 
The $PQ$-relaxation can be strengthned by appending valid constraints derived with the reformulation linearization (RLT) technique \citep{Sherali1992}.

{\footnotesize
\begin{subequations}
\label{Eq:PQ_Model}
\begin{align}
\min_{\vect{y}, \vect{z}, \vect{v}, \vect{q}}
& \sum_{(i,l,j)\in I\times L\times J} (\ci - d_j) \vilj + \sum\limits_{(i,j)\in Z} (\ci - d_j) \zij \label{Eq:PQ_Model_Objective} \\
& \Ail \leq \sum_{(l,j)\in Y} \vilj + \sum\limits_{j \in J} \zij \leq  \Aiu & i \in I \label{Eq:PQ_Model_Supply_Bounds} \\
& \sum_{(i,j) \in Z} \vilj \leq \Sl & l \in L \label{Eq:PQ_Model_Pool_Capacity} \\
& \Djl \leq \sum_{(i,l)\in X} \vilj + \sum_{i \in I} \zij \leq \Dju & j \in J \label{Eq:PQ_Model_Demand_Bounds} \\
& \vilj = \qil \ylj & i \in I, j \in J, l \in L \label{Eq:PQ_Model_Path_Flow} \\
& \sum_{i \in I} \qil = 1 & l \in L \label{Eq:PQ_Model_Total_Fractions} \\
& \sum_{i \in I} \vilj = \ylj & l\in L, j\in Y \label{Eq:PQ_Model_Output_Flow} \\
& \sum_{(i, l)\in X} ( \Cik - \Pjkl ) \vilj + \sum_{i \in I} ( \Cik - \Pjkl ) \zij \geq 0 & j\in J \label{Eq:PQ_Model_Quality_Lower} \\
& \sum_{(i, l)\in X} ( \Cik - \Pjku ) \vilj + \sum_{i \in I} ( \Cik - \Pjku ) \zij \leq 0 & j\in J \label{Eq:PQ_Model_Quality_Upper} \\
& \ylj, \zij \geq 0 & i\in I, l\in L, j\in J \label{Eq:PQ_Model_Flow_Bounds} \\
& 0 \leq \qil \leq 1 & i \in I, l\in L \label{Eq:PQ_Model_Fraction_Bounds} 
\end{align}
\end{subequations}
}
%
Expression (\ref{Eq:PQ_Model_Objective}) minimizes the total cost. 
Constraints (\ref{Eq:PQ_Model_Supply_Bounds}) - (\ref{Eq:PQ_Model_Demand_Bounds}) enforce bounds on the raw material quantities, pool sizes, and final product quantities, respectively.
Constraints (\ref{Eq:PQ_Model_Path_Flow}) - (\ref{Eq:PQ_Model_Output_Flow}) express material balances.
Constraints (\ref{Eq:PQ_Model_Quality_Lower}) - (\ref{Eq:PQ_Model_Quality_Upper}) impose quality specifications.
Constraints (\ref{Eq:PQ_Model_Flow_Bounds}) - (\ref{Eq:PQ_Model_Fraction_Bounds}) ensure non-negativity of flow variables and fraction bounds.

\subsection{Computational Complexity and Approximation Algorithms} 
\label{ss:pp_complexity}

This section discusses the known computational complexity and approximation algorithms for the pooling problem. Table \ref{tab:pooling_complexity} additionally summarizes the computational complexity results discussed in this section.

When there are no quality constraints, i.e.\ $|K|=0$, or $P_{j,k}^L=0$ and $P_{j,k}^U=+\infty$ for each $j\in J$ and $k\in K$, pooling  becomes an instance of the well-known minimum cost flow problem which is polynomially solvable.
Pooling is also a tractable LP when there are no intermediate pools and the problem is referred to as blending.
In the more general case with both quality constraints and intermediate pools, Table~\ref{tab:pooling_complexity} reports state-of-the-art computational complexity results for subproblems with (i) set cardinality restrictions, (ii) special network structure and (iii) supply/demand/capacity restrictions.

\citet{Alfaki2013formulations} show that pooling is strongly $\mathcal{NP}$-hard even in the special case with a single pool, i.e.\ $|L|=1$, through a reduction from the independent set problem.
The problem remains $\mathcal{NP}$-hard for instances with a single quality attribute, i.e.\ $|K|=1$, via a reduction from Exact Cover by 3-Sets \citep{Boland2017}.
On the other hand, in the singleton cases with a single input or output, i.e.\ $\min\{|I|,|J|\}=1$, pooling can be easily formulated as an LP and is therefore polynomially solvable \citep{Dey2015}.
These findings have motivated further investigations on pooling with set cardinality restrictions.
For instances with a single pool with no input-output arcs where the number of inputs \citep{Boland2017}, outputs \citep{Alfaki2013formulations}, or attributes \citep{Alfaki2013formulations, Haugland2014} is bounded by a constant, i.e.\ $\min\{|I|,|J|,|K|\}=O(1)$, the problem can be solved in polynomial time using a series of LPs.
In the case $|K|=1$ with a single quality attribute where there are either two inputs, or two outputs, i.e.\ $\min\{|I|,|J|\}=2$, the problem is still $\mathcal{NP}$-hard by a reduction from Exact Cover by 3-Sets.
Finally, when $|I|=|J|=2$, pooling is known to be weakly $\mathcal{NP}$-hard through a reduction from Partition.

\citet{Haugland2016} shows that pooling is $\mathcal{NP}$-hard for problem instances with sparse network structure.
Let $\Delta_i^{\stout}$ and $\Delta_l^{\stout}$ be the out-degree, i.e.\ number of outgoing arcs, of input $i\in I$ and pool $l\in L$, respectively.
When every out-degree is at most two, i.e.\ $\max\{\Delta_i^{\stout}, \Delta_l^{\stout}\} \leq 2$, \citet{Haugland2016} presents an $\mathcal{NP}$-hardness reduction from maximum satisfiability.
Denote by $\Delta_l^{\stin}$ and $\Delta_j^{\stin}$ the in-degree, i.e.\ number of ingoing arcs, of pool $l\in L$ and output $j\in J$, respectively.
When each in-degree does not exceed two, i.e.\ $\max\{\Delta_l^{\stin}, \Delta_j^{\stin}\} \leq 2$, pooling is $\mathcal{NP}$-hard through a reduction from minimum satisfiability \citep{Haugland2016}.
However, in the case where each pool has either in-degree or out-degree equal to one, i.e.\ $\min\{\Delta_l^{\stout},\Delta_l^{\stin}\}=1$, the problem is polynomially solvable \citep{Dey2015, Haugland2016polynomial}. Finally, for instances with a single pool and attribute, unlimited supplies/pool capacities and fixed demands, the pooling problem is strongly-polynomially solvable \citep{baltean2018piecewise}.

\begin{table}
\footnotesize
\caption{Pooling Problem Computational Complexity Results}
\begin{center}
\begin{tabular}{|l c c|}
\hline
\bf Subproblem & \bf Complexity & \bf Reduction \\
\hline
\multicolumn{3}{|l|}{\bf Singleton subproblems}\\
$|I|=1$ & $\mathcal{P}$ & - \\
$|J|=1$ & $\mathcal{P}$ & - \\
$|L|=1,\ Z=\emptyset$ & $\mathcal{NP}$-hard & Independent Set \\
$|K|=1$ & $\mathcal{NP}$-hard & Exact Cover by 3-Sets \\
\hline
\multicolumn{3}{|l|}{\bf Other cardinality-restricted special cases}\\
$|L|=1,\ |I|=O(1),\ Z=\emptyset$ & $\mathcal{P}$ & - \\
$|L|=1,\ |J|=O(1),\ Z=\emptyset$ & $\mathcal{P}$ & - \\
$|L|=1,\ |K|=O(1),\ Z=\emptyset$ & $\mathcal{P}$ & - \\
$|I|=2,\ |K|=1$ & $\mathcal{NP}$-hard & Exact Cover by 3-Sets \\
$|J|=2,\ |K|=1$ & $\mathcal{NP}$-hard & Exact Cover by 3-Sets \\
$|I|=2,\ |J|=2,\ |K|=1$ & $\mathcal{NP}$-hard & Partition \\
\hline
\multicolumn{3}{|l|}{\bf Special network structure}\\
$\min\{\Delta_l^{\stout},\ \Delta_l^{\stin}\} = 1$ & $\mathcal{P}$ & - \\
$\Delta_i^{\stout}\leq 2,\ \Delta_l^{\stout}\leq 2$ & $\mathcal{NP}$-hard & Maximum Satisfiability \\
$\Delta_l^{\stin}\leq 2,\ \Delta_j^{\stin}\leq 2$ & $\mathcal{NP}$-hard & Minimum Satisfiability \\
\hline
\multicolumn{3}{|l|}{\bf Supply, demand, and pool capacity restrictions}\\
$|L|=1,\ |K|=1,\ A^L_i = 0, A^U_i = \infty,\ S_l = \infty,\ D^L_j=D^U_j$ & $\mathcal{P}$ & - 
\\ 
\hline
\end{tabular}
\end{center}
\label{tab:pooling_complexity}
\end{table}



The only known theoretical performance bounds for pooling are an $O(n)$-approximation algorithm and an $\Omega(n^{1-\epsilon})$ inapproximability result, for $\epsilon>0$, by \citet{Dey2015}.
The proposed algorithm solves the relaxation obtained by applying piecewise linear McCormick envelopes to the bilinear terms of the $PQ$-formulation \citep{Gupte2013, Gupte2017}.
Overestimators and underestimators are computed by partitioning the domain of proportion variables to a finite MILP-representable set. 
For the negative result, \citet{Dey2015} present an approximation-preserving reduction from independent set.


\section{Process Scheduling}
\label{Section:Process_Scheduling}

\begin{figure}
    \centering
    \includegraphics[scale=0.85]{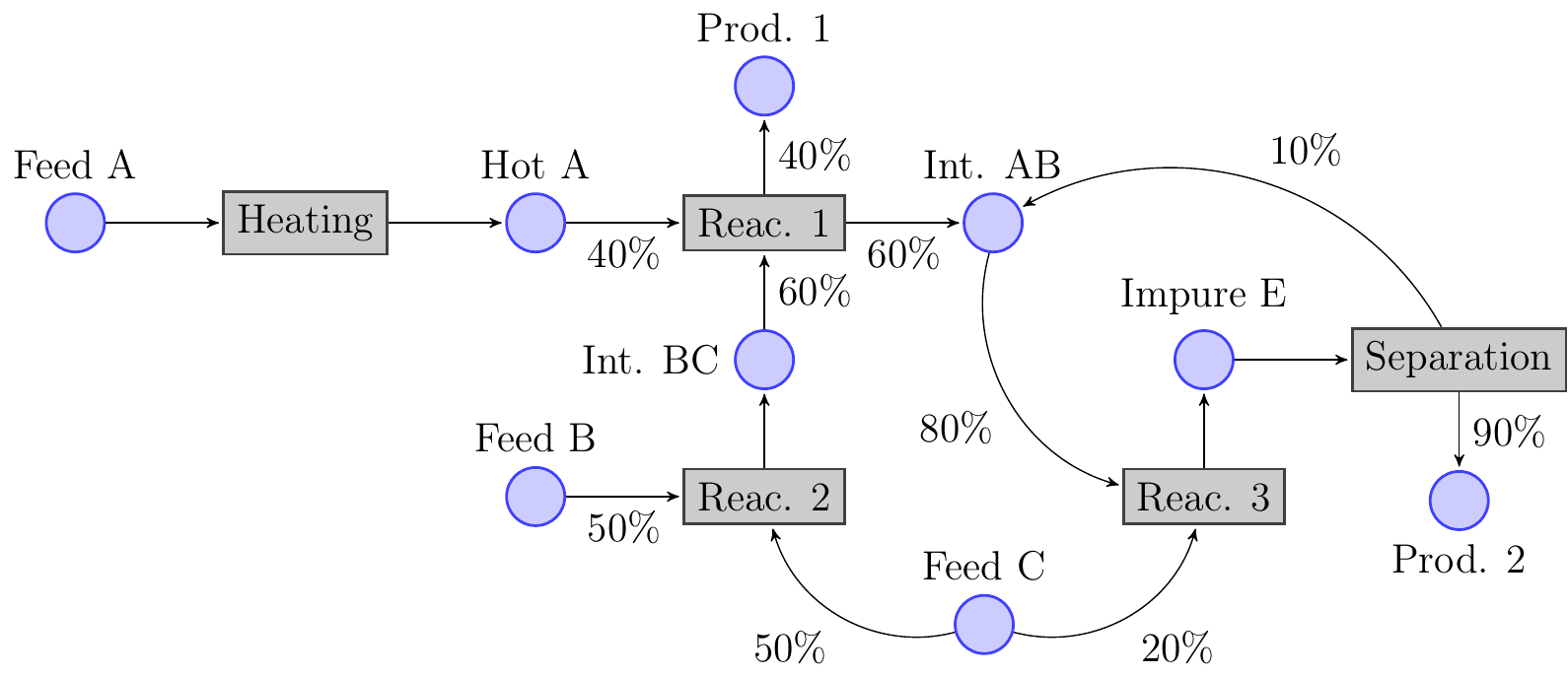}
    \caption{Network structure for the state-task network by Kondili et al.\citep{Kondili1993}}
    \label{fig:stn}
\end{figure}

Scheduling process operations, a.k.a.\ batch scheduling, is crucial in different application areas including chemical manufacturing, pharmaceutical production \citep{Lainez2012}, food industry \citep{Stefanis1997}, and oil refining \citep{Harjunkoski2014}. 
Process scheduling problems are the topic of many fruitful investigations in the PSE literature \citep{Castro2018, Floudas2004, Harjunkoski2014, Mendez2006, Wiebe2018}. 
The goal is to efficiently allocate the limited resources, e.g.\ processing units, of manufacturing plants to tasks and decide the product batch sizes so as to construct multiple intermediate and final products satisfying the customer demand.
These products are often based on recipes in the form of \emph{state-task networks} where each task receives raw materials and intermediate products to generate new products \citep{Kondili1993, Shah1993}.
State-task networks may model general batch processes including material mixing, splitting, recycling, as well as different storage policies \citep{Kallrath2002}.
Fig.~\ref{fig:stn} presents an example of a state-task network. 
Typically, process scheduling involves solving $\mathcal{NP}$-hard, mixed-integer linear programming problems which require algorithms exploiting the state-task network's structure. 



\subsection{Brief Literature Overview}

Scheduling is a relatively recent area in PSE \citep{Mauderli1979, Reklaitis1982} and has received considerable attention after the seminal work by \citet{Kondili1993} who introduced the \emph{state-task network} framework for modeling mixing and splitting of material batches.
\citet{Pantelides1994} extended the state-task network to the notion of a Resource-Task Network (RTN) for incorporating multiple resources in a unified setting.
Process scheduling problems include a variety of aspects that need to be considered, such as different production stages, storage policies, demand patterns, changeovers, resource constraints, time constraints, and uncertainty.
Furthermore, they require optimizing different objective functions, e.g.\ makespan, production costs, or sales profit.
There is significant work providing surveys and problem classification for process scheduling \citep{Castro2018, Floudas2004, Harjunkoski2014, Li2008, Maravelias2012, Mendez2006}.

Significant literature solves process scheduling problems using MILP, this work is supported by the significant progress in CPU speed and algorithms in the last two decades.
State-of-the-art mathematical modeling develops discrete-time and continuous-time formulations \citep{Maravelias2003continuous, Maravelias2003discrete}.
These approaches are strengthened by reformulation and tightening methods \citep{Ierapetritou1998, Schilling1996, Sundaramoorthy2005, Velez2013}.
Branch-and-cut, decomposition, constraint programming, metaheuristic, hybrid approaches, and satisfiability modulo theories are also explored \citep{Castro2011,kopanos2009bi,10.1007/978-3-540-24664-0_1,MISTRY201898,TILL2007630,VELEZ201328,wu2003decomposition}.
Recently, generalized-disjunctive programming has emerged as a novel framework for effectively solving process scheduling problems using big-M and convex hull reformulations \citep{Castro2012}.
In addition, rescheduling has been used as a tool for mitigating the effect of disturbances under uncertainty \citep{Gupta2016,Gupta2019}.

\subsection{Problem Definition}

A process scheduling problem instance consists of a state-task network specifying a recipe for generating chemical products from raw materials.
Formally, a state-task network is a directed bipartite graph $N=(S\cup I,A)$ with a partition of nodes into a subset $S$ corresponding to states, i.e.\ raw materials, intermediate, and final products, and a subset $I$ representing tasks.
The network $N$ is bipartite, i.e.\ the set $A=A^-\cup A^+$ of arcs consists of \emph{consumption arcs} $A^-\subseteq S\times I$ and \emph{production arcs} $A^+\subseteq I\times S$.
Arc $(s,i)\in A^-$ implies that task $i\in I$ consumes a positive amount of state $s\in S$. 
Analogously, arc $(i,s)\in A^+$ indicates that task $i\in I$ produces a positive amount of state $s\in S$.  
There is a set $J$ of processing units for executing tasks.
Each unit may process at most one task per unit of time.
Denote by $J_i\subseteq J$ the subset of units that may perform task $i \in I$ and by $I_j=\{i:(j\in J_i)\wedge(i\in I)\}$ the set of tasks that may be performed by unit $j\in J$.
Moreover, let $\ell$ be the time horizon length and $T=[0,\ell]$.
If task $i\in I$ begins execution on unit $j\in J_i$ at time $t\in T$, then it may process a variable amount $b_{i,j,t}$ of material, a.k.a.\ \emph{batch size}, for $p_{i,j}$ units of time, and completes at time $t+p_{i,j}$.
Let $b_{i,j}^L$ and $b_{i,j}^U$ be the minimum and maximum capacity, respectively, of unit $j\in J$ when processing task $i\in I$.
Continuous variable $b_{i,j,t}$ is allowed to take any value in the interval $[b_{i,j}^L, b_{i,j}^U]$.
Denote by $f_{i,s}^-$ and $f_{i,s}^+$ the material fraction entering and exiting processing, respectively, as state $s\in S$ when task $i\in I$ is processed.
If $S_i^-=\{s: ((s,i)\in A^-) \wedge (s\in S) \}$ and $S_i^+=\{s: ((i,s)\in A^+) \wedge (s\in S)\}$ are the consumables and products of task $i\in I$, respectively, then task $i$ consumes $f_{i,s}^-b_{i,j,t}$ portion of state $s\in S_i^-$ and produces $f_{i,s}^+b_{i,j,t}$ quantity of state $s\in S_i^+$.
Suppose that $I_s^+=\{i: ((i,s)\in A^+)\wedge(i\in I)\}$ and $I_s^-=\{i: ((s,i)\in A^-)\wedge(i\in I)\}$ are the sets of tasks producing and consuming, respectively, state $s\in S$.
Then, $\sum_{i\in I_s^+}\sum_{j\in J_i}\sum_{t\in T} f_{i,s}^+b_{i,j,t-p_{i,j}}$ amount is produced and 
$\sum_{i\in I_s^-}\sum_{j\in J_i}\sum_{t\in T} f_{i,s}^-b_{i,j,t}$ quantity is consumed for state $s\in S$ at $t\in T$.
\ref{app:scheduling_notation} presents the notation for process scheduling.


The goal of process scheduling is to satisfy a demand $d_s$ for each state $s\in S$.
Denote by $y_{s,t}$ the amount of $s\in S$ at time slot $t\in T$.
Without loss of generality, we assume that $y_{s,0}=0$, i.e.\ there is initially zero amount of state $s\in S$.
When the time horizon completes, the obtained solution must satisfy $y_{s,\ell}\geq d_s$ for each $s\in S$.
The objective is to schedule the tasks on the units and decide the batch sizes so that the makespan $z$, i.e.\ the time at which the last task completes, is minimized. 


\subsection{Mathematical Models}

The main approaches for formulating process scheduling problems as MILP problems are typically classified as (i) discrete-time \citep{Kondili1993,Shah1993}, or (ii) continuous-time \citep{Maravelias2003continuous}.
\citet{Floudas2004} and \citet{Mendez2006} provide thorough discussions on the advantages of each. 
Discrete-time formulations partition time into a large number of time intervals.
Continuous-time formulations (i) use a small number of event points 
resulting in fewer variables, and (ii) express inventory and backlog costs linearly.
However, continuous-time formulations generally tend to be nonlinear. 
Mixed-time representations utilize both the discrete-time and continuous-time models \citep{Lee2018,Maravelias2005}.


\subsubsection{Discrete-Time Formulation}

Discrete-time formulations partition the time horizon into a set $T=\{1,\ldots,\ell\}$ of equal-length slots.
Integer variable $x_{i,j,t}$ indicates whether task $i\in I$ is executed by unit $j\in J$ starting at time $t\in T$.
Continuous variable $b_{i,j,t}$ specifies the corresponding batch size.
Continuous variables $y_{s,t}$ denote the stored amount of state $s\in S$ at time $t\in T$. 
Finally, continuous variable $z$ computes the makespan.
Process scheduling can be modeled using the Eq.\ (\ref{Eq:DT_Model}) MILP formulation.

{\allowdisplaybreaks
{\footnotesize
\begin{subequations}
\label{Eq:DT_Model}
\begin{align}
\min \quad & z
\label{Eq:DT_Model_Objective} \\
& z \geq x_{i,j,t} (t + p_{i,j}) & i\in I, j\in J_i, t\in T \label{Eq:DT_Model_Makespan} \\
& \sum_{i \in I_j} \sum_{t' = t-\pij+1}^{t} x_{i,j,t'} \leq 1 & j \in J, t \in T \label{Eq:DT_Model_Task_Capacity} \\
& x_{i,j,t} b_{i,j}^L \leq b_{i,j,t} \leq x_{i,j,t} b_{i,j}^U & i \in I, j \in J_i, t \in T \label{Eq:DT_Model_Material_Capacity} \\
& y_{s,t} = y_{s,t-1} + \sum_{i \in I^+_s} \sum_{j \in J_i} f_{i,s}^+ b_{i,j,t-\pij+1} - 
\sum_{i \in I_s^-} \sum_{j \in J_i} f_{i,s}^- b_{i,j,t} & s \in S, t \in T \label{Eq:DT_Model_Material_Balance} \\
& y_{s,\ell} \geq d_s, y_{s,1}=0 & s\in S \label{Eq:DT_Model_State_Capacity} \\
& x_{i,j,t} \in \{0,1\} & i\in I, j\in J_i, t\in T \label{Eq:DT_Model_Integrality} \\
& z, b_{i,j,t}, y_{s,t} \geq 0 & i\in I, j\in J_i, s \in S, t \in T \label{Eq:DT_Model_Positiveness}
\end{align}
\end{subequations}
}}
Expression (\ref{Eq:DT_Model_Objective}) minimizes makespan.
Constraints (\ref{Eq:DT_Model_Makespan}) define the makespan.
Constraints (\ref{Eq:DT_Model_Task_Capacity}) ensure that each unit processes at most one task at each point in time.
Constraints (\ref{Eq:DT_Model_Material_Capacity}) and (\ref{Eq:DT_Model_Material_Balance}) express unit capacities and material conservation, respectively.
Constraints (\ref{Eq:DT_Model_State_Capacity}) enforce that the demand is statisfied.
Finally, constraints (\ref{Eq:DT_Model_Integrality}) - (\ref{Eq:DT_Model_Positiveness}) impose that integer and continuous variables are binary and non-negative, respectively.

\subsubsection{Continuous-Time Formulation}

Continuous-time formulations divide the time horizon into a set of slots, similarly to discrete-time formulations.
The number of slots is fixed, but the slots are not necessarily of equal length. 
The slot boundaries are determined by a set $T$ of $\ell$ variable time points. 
Continuous variable $t_k$ specifies the rightmost time point $k\in T$ of one slot and the leftmost time point of the subsequent slot. 
Furthermore, each job's starting and completion time is mapped to a time point.
Binary variables $x_{i,k}^S$ and $x_{i,k}^F$ express whether task $i\in I$ begins and completes, respectively, at time point $k\in T$.
To match time points with task starting and completion times, continuous variables $t_{i,k}^S$ and $t_{i,k}^F$ compute the start and finish time of task $i\in I$ beginning at time point $k\in T$.
Continuous variables $p_{i,k}$ and $b_{i,k}$ correspond to the processing time and batch size of task $i\in I$ starting at time point $k\in T$. 
Finally, continuous variable $y_{s,k}$ models the amount of state $s\in S$ at time point $k\in T$. 
Without loss of generality, we assume that $|J_i|=1$ for each $i\in I$, i.e.\ tasks $i$ can only be executed by a single unit.
To model the case $|J_i|>1$, we may add multiple occurrences of the same task.
Furthermore, we note that continuous-time formulations may easily incorporate variable task durations.
Specifically, we suppose that task $i\in I$ has a variable duration $\alpha_i$ that depends on the batch size, in addition to a fixed duration $\beta_i$.
Then, variable $p_{i,k}$ denotes the processing time of job $i\in I$ starting at time point $k\in T$. 

{\allowdisplaybreaks
{\footnotesize
\begin{subequations}
\label{Eq:CT_Model}
\begin{align}
\min \quad & z \label{Eq:CT_Model_Objective} \\
& z \geq t_{i,k}^S + p_{i,k} & i\in I, k\in T \label{Eq:CT_Model_Makespan} \\
& p_{i,k} = \alpha_i x_{i,k}^S + \beta_i b_{i,k}^S & i\in I, k\in T \label{Eq:CT_Model_Processing_Time} \\
& t_{i,k}^S \leq t_k + H(1-x_{i,k}^S) & i\in I, k\in T \label{Eq:CT_Model_Binary_Activation_1} \\
& t_{i,k}^S \geq t_k - H(1-x_{i,k}^S) & i\in I, k\in T \label{Eq:CT_Model_Binary_Activation_2} \\
& t_{i,k}^F \leq t_k + p_{i,k} + H(1-x_{i,k}^S) & i\in I, k\in T \label{Eq:CT_Model_Binary_Activation_3} \\
& t_{i,k}^F \geq t_k + p_{i,k} - H(1-x_{i,k}^S) & i\in I, k\in T \label{Eq:CT_Model_Binary_Activation_4} \\
& t_{i,k}^F - t_{i,k-1}^F \leq H x_{i,k}^S & i\in I, k\in T \setminus \{1\}\\
& t_{i,k-1}^F \leq t_{k} + H(1-x_{i,k}^F) & i\in I, k\in T \setminus \{1\} \label{Eq:CT_Model_Binary_Activation_5} \\
& t_{i,k-1}^F \geq t_{k} - H(1-x_{i,k}^F) & i\in I, k\in T \setminus \{1\} \label{Eq:CT_Model_Binary_Activation_6} \\
& t_1=0, t_{k-1} \leq t_{k}, t_{\ell} = H & k\in T \setminus \{1\} \label{Eq:CT_Model_Increasing_Times} \\
& \sum_{i\in I_j} \sum_{k' \leq k} (x_{i,k'}^S - x_{i,k}^F) \leq 1 & j\in J, k\in T \label{Eq:CT_Model_Task_Capacity} \\
& \sum_{k\in T} x_{i,k}^S = \sum_{k\in T} x_{i,k}^F & i\in I \label{Eq:CT_Model_Task_Completion} \\
& x_{i,k}^S b_i^L \leq b_{i,k} \leq x_{i,k}^S b_i^U & i\in I, k\in T \label{Eq:CT_Model_Unit_Capacity_1} \\
& y_{s,k} = y_{s,k-1} + \sum_{i\in I_s^+} f_{i,s}^+b_{i,k-1} - \sum_{i\in I_s^-} f_{i,s}^- b_{i,k} & s\in S, k\in T\setminus\{1\} \label{Eq:CT_Model_Mass_Balance_2}  \\
& y_{s,\ell} \geq d_s &  s\in S, k\in T \label{Eq:CT_Model_Storage_Capacity}  \\
& x_{i,k}^S,x_{i,k}^F\in\{0,1\} & i\in I, k\in T \label{Eq:CT_Model_Integrality} \\
& t_k, t_{i,k}^S, t_{i,k}^F, p_{i,k}, b_{i,k}, y_s \geq 0 & i\in I, k\in T, s\in S \label{Eq:CT_Model_Positiveness} 
\end{align}
\end{subequations}
}}
Expression (\ref{Eq:CT_Model_Objective}) minimizes makespan.
Constraints (\ref{Eq:CT_Model_Makespan}) are the makespan definition.
Constraints (\ref{Eq:CT_Model_Processing_Time}) calculate the job processing time.
Binary activation constraints (\ref{Eq:CT_Model_Binary_Activation_1}) - (\ref{Eq:CT_Model_Binary_Activation_6}) map continuous time variables to time points.
Constraints (\ref{Eq:CT_Model_Increasing_Times}) impose time horizon boundaries and the non-decreasing order of time points.
Constraints (\ref{Eq:CT_Model_Task_Capacity}) enforce that each unit processes at most one task per time.
Constraints (\ref{Eq:CT_Model_Task_Completion}) ensure that every task that begins processing must also complete. 
Constraints (\ref{Eq:CT_Model_Unit_Capacity_1}) incorporate unit capacities.
Constraints (\ref{Eq:CT_Model_Mass_Balance_2}) express mass balance.
Constraints (\ref{Eq:CT_Model_Storage_Capacity}) model storage capacities.
Finally, constraints (\ref{Eq:CT_Model_Integrality}) - (\ref{Eq:CT_Model_Positiveness}) ensure that continuous and integer variables are non-negative and binary, respectively.

\subsection{Computational Complexity and Approximation Algorithms}

Process scheduling problems are frequently characterized as computationally challenging \citep{Floudas2004, Harjunkoski2014}.
However, computational complexity investigations are limited and isolated.
To our knowledge, \citet{Burkard1998} have only work in this direction. \citet{Burkard1998} observe that process scheduling (i) is strongly $\mathcal{NP}$-hard as a generalization of the job shop scheduling problem, and (ii) remains $\mathcal{NP}$-hard even in the special case with two states through a straightforward reduction from knapsack.
Heuristics have been reported as a tool for solving large-scale process scheduling instances \citep{Harjunkoski2014, Mendez2006, Panwalkar1977}.
Nevertheless, only few early works in the area develop heuristics exploiting the problem's combinatorial structure, e.g.\ greedy layered \citep{Blomer2000} and discrete-time relaxation rounding \citep{Burkard1998}. 
Furthermore, there is lack of analytically proven performance guarantees.



The above observations are opposed to the tremendous contributions of computational complexity and approximation algorithms in scheduling theory.
A classical scheduling problem may be defined using the \emph{three-field notation} which incorporates \citep{Graham1979}:
(i) a machine environment, (ii) job characteristics, and (iii) an objective function.
The goal is to decide when and where to execute the jobs, i.e.\ at which times and on which machines, so that the objective function is optimized. 
Single stage machine environments include: identical, related, and unrelated machines. 
Multistage machine environments can be: open shops, flow shops, or job shops.
Examples of job characteristics are: release times, deadlines, and precedence constraints.
Objective functions include: makespan, response time, tardiness, throughput and others.
Despite the commonalities between process scheduling and classical scheduling theory, there is a striking absence of connections between the two fields.
Their synergy constitutes a particularly interesting future direction and has strong potential for successfully solving open process scheduling problems.
To this end, \citet{Chen1998} and \citet{SchedulingZoo} provide an extensive survey of results for fundamental scheduling problems.
\citet{Brucker2006, Leung2004} and \citet{Pinedo2012} present a collection of algorithms and techniques for effectively solving such problems.

\section{Heat Exchanger Network Synthesis}
\label{Section:Heat_Exchangers}
Heat exchanger network synthesis is one of the most extensively studied problems in chemical engineering \citep{Biegler1997, Escobar2013, Furman2002, Gundersen1988, Smith2000}. 
Major heat exchanger network synthesis applications include energy systems producing liquid transportation fuels \citep{Floudas2012, Niziolek2015}, natural gas refineries \citep{Baliban2010, Mehdizadeh2017}, refrigeration systems \citep{Shelton1986}, batch processes \citep{Castro2015, Zhao1998}, and water utilization systems \citep{Bagajewicz2002}.
Heat exchanger network synthesis minimizes the total investment and operating costs in chemical processes.
In particular, heat exchanger network synthesis: (i) improves energy efficiency by reducing heating utility usage, (ii) optimizes network costs by accounting for the number of heat exchanger units and area physical constraints, and (iii) improves energy recovery by integrating hot and cold process streams \citep{Elia2010, Floudas1987}.
The goal is designing a heat exchanger network matching hot streams to cold streams and recycling residual heat, by taking into account the nonlinear nature of heat exchange and thermodynamic constraints. 
Heat exchanger network synthesis is an $\mathcal{NP}$-hard, MINLP instance with (i) nonconvex nonlinearities for enforcing energy balances, and (ii) discrete decisions for placing heat exchanger units.
This section investigates the nonlinear and integer heat exchanger network synthesis parts individually by considering the multistage minimum utility cost, and minimum number of matches problems separately.

\subsection{Brief Literature Overview}

Optimization methods for heat exchanger network synthesis can be classified as: (i) simultaneous, or (ii) sequential.
Simultaneous methods produce globally optimal solutions. 
Sequential methods do not provide any guarantee of optimality, but are useful in practice.
Simultaneous methods formulate heat exchanger network synthesis as a single MINLP, e.g.\ \citet{Papalexandri1994}.
\citet{Ciric1991} propose the \emph{hyperstructure MINLP} formulating heat exchanger network synthesis without decomposition based on the stream superstructure introduced by \citet{Floudas1986}.
\citet{Yee1990} develop the \emph{multistage MINLP} (a.k.a.\ SYNHEAT model) using a stagewise superstructure.
Because the multistage MINLP assumes isothermal mixing at each stage, the nonlinear heat balances are simplified and performed only between stages.
Sequential methods decompose heat exchanger network synthesis into three distinct subproblems: (i) minimum utility cost, (ii) minimum number of matches, and (iii) minimum investment cost.
These subproblems are more tractable than simultaneous heat exchanger network synthesis.
In particular, \citet{Cerda1983matches, Cerda1983utility, Papoulias1983} suggest the transportation and transshipment models formulating the minimum utility cost problem as LP and the minimum number of matches problem as MILP. 
\citet{Floudas1986} propose the stream superstructure formulating the minimum investment cost problem as an NLP.
Other heat exchanger network synthesis approaches exploit the problem's thermodynamic nature, and mathematical and physical insights in order to design more efficient algorithms. 
\citep{Ahmad1989supertargeting, Ahmad1989shells, Gundersen1990, Gundersen1997, Kouyialis2017, Leitold2019, Linnhoff1989, Linnhoff1978, Linnhoff1983, Masso1969, Mistry20161, Pho1973, Polley1999}.



\subsection{Problem Definitions}

A heat exchanger network synthesis instance consists of a set $H$ of hot streams to be cooled down and a set $C$ of cold streams to be heated up.
Each hot stream $i\in H$ and cold stream $j\in C$ is associated with an initial, inlet temperature $T_i^{\stin}$, $T_j^{\stin}$, target, outlet temperature $T_i^{\stout}$, $T_j^{\stout}$, and flow rate heat capacity $F_i$, $F_j$, respectively.
The temperature of hot stream $i\in H$ must be decreased from $T_i^{\stin}$ down to $T_i^{\stout}$, while the temperature of cold stream $j\in C$ has to be increased from $T_j^{\stin}$ up to $T_j^{\stout}$.
For each $i\in H$ and $j\in C$, flow rate heat capacities $F_i$ and $F_j$ specify the quantity of heat that a stream releases and absorbs, respectively, per unit of temperature change.
That is, hot stream $i\in H$ supplies $F_i(T_i^{\stin}-T_i^{\stout})$ units of heat, while cold stream $j\in C$ demands $F_j(T_j^{\stout}-T_j^{\stout})$ units of heat.
\ref{app:hens_notation} presents the notation for heat exchanger network synthesis.

\subsubsection{Multistage Minimum Utility Cost}

\begin{figure*}[t!]
\centering
\includegraphics[]{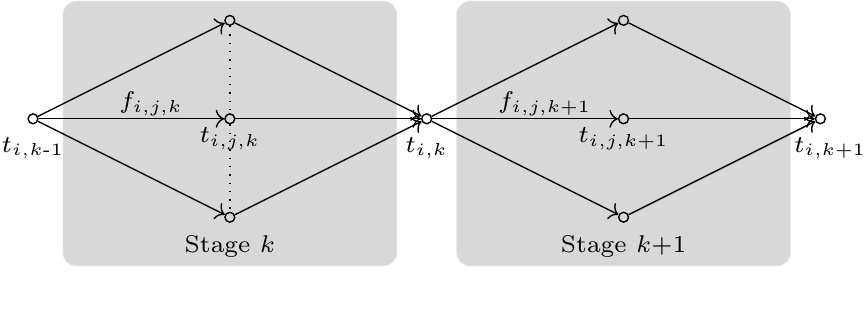}
\caption{
Illustration of multistage heat exchanger network synthesis. 
Hot stream $i\in H$ across multiple stages in increasing order of their indices.
At stage $k\in S$, stream $i$ is split into substreams.
}
\label{Fig:multistage}
\end{figure*}

In multistage heat exchanger network synthesis, heat transfers between streams occur in a set $S$ of $\ell$ different stages.
Hot streams flow from the stage $1$ to stage $\ell$, while cold streams flow, in the opposite direction, from stage $\ell$ to stage $1$. 
When a hot, respectively cold, stream enters stage $k\in S$, it is split into substreams each one exchanging heat with exactly one cold, respectively hot, stream and these substreams are merged back together when the stream exits the stage.
Figure~\ref{Fig:multistage} illustrates splitting and mixing.
For $k\in S$, denote by $t_{i,k}$ the temperature of hot stream $i\in H$ when exiting and entering the stages $k$ and $k+1$, respectively.
Similarly, let $t_{j,k}$ be the initial and last temperature of cold stream $j\in H$ at stages $k$ and $k+1$, respectively, for $k\in S$.
The multistage minimum utility cost problem decides how to split the streams in each stage.
The substream of $i\in H$ exchanging heat with $j\in C$ at stage $k\in S$ gets flow rate heat capacity $f_{i,j,k}^H$.
Similarly, the substream of $j\in C$ exchanging heat with $i\in H$ at stage $k\in S$ is assigned flow rate heat capacity $f_{i,j,k}^C$.
It must be the case that $\sum_{j\in C}f_{i,j,k}^H=F_i$ and $\sum_{i\in H}f_{i,j,k}^C=F_j$, for all $k\in S$.
If $i\in H$ is matched with $j\in C$ at $k\in S$, the corresponding substream of $i$ and $j$ results with a temperature $t_{i,j,k}^H$ and $t_{i,j,k}^C$, respectively, when the stage completes.
At stage $k\in S$, hot stream $i\in H$ and cold stream $j\in C$ have final temperatures $t_{i,k}$ and $t_{j,k-1}$ such that 
$F_it_{i,k}=\sum_{j\in C}f_{i,j,k}^Ht_{i,j,k}^H$ and $F_jt_{j,k-1}=\sum_{i\in H}f_{i,j,k}^Ct_{i,j,k}^C$.
The total heat exchanged between $i$ and $j$ at $k$ is $q_{i,j,k}=f_{i,j,k}^H(t_{i,k}-t_{i,j,k}^H)$ and $q_{i,j,k}=f_{i,j,k}^C(t_{i,j,k}^C-t_{j,k})$, i.e.\ there is heat conservation.
A hot and cold utility may provide or extract heat at unitary costs $c^{HU}$ and $c^{CU}$.
The cold utility exports $Q_i^{CU} = F_i(t_{i,\ell} - T_i^{\text{out}})$ units of heat from hot stream $i\in H$.
Analogously, the hot utility supplies $Q_j^{HU} = F_j(T_j^{\text{out}} - t_{j,0})$ units of heat to cold stream $j\in C$.
The goal is to exchange heat and reach the target temperature for each stream so that the total utility cost $\sum_{i\in H}c^{CU}Q_i^{CU} + \sum_{j\in C}c^{HU}Q_j^{HU}$ is minimized.

\subsubsection{Minimum Number of Matches Problem}

In the minimum number of matches problem, heat transfers occur similarly to standard network flow problems \citep{Ahuja1993}. 
A problem instance only consists of streams. 
The utilities are considered as streams whose parameters, i.e.\ flow rate heat capacities, inlet and outlet temperatures, are computed by solving a minimum utility cost LP to ensure heat conservation.
Specifically, hot stream $i\in H$ exports $h_i=F_i(T_i^{\text{in}}-T_i^{\text{out}})$ units of heat, cold stream $j\in C$ receives $c_j=F_j(T_j^{\text{out}}-T_j^{\text{in}})$ units of heat, and $\sum_{i\in H}h_i=\sum_{j\in C}c_j$.
A minimum heat approach temperature $\Delta T_{\min}$ accounts for the energy lost by the system.
We may assume that $\Delta T_{\min}=0$, because any problem instance can be transformed to an equivalent one satisfying this assumption.
Let $T_0>T_1>\dots>T_r$ be all discrete inlet and outlet temperature values. 
The temperature range is partitioned into a set $T=\{[T_t,T_{t-1}]:1\leq t\leq r\}$ of consecutive temperature intervals.
In temperature interval $t\in T$, hot stream $i\in H$ exports $\sigma_{i,t}=F_i(T_{t-1}-T_t)$ units if 
$[T_t,T_{t-1}]\subseteq[T_i^{\text{out}},T_i^{\text{in}}]$, and $\sigma_{i,t}=0$ otherwise. 
Likewise, cold stream $j\in C$ receives $\delta_{j,t}=F_j(T_{t-1}-T_t)$ units of heat if $[T_t,T_{t-1}]\subseteq[T_j^{\text{in}},T_j^{\text{out}}]$, and $\delta_{j,t}=0$ otherwise.
A feasible solution specifies a way to transfer the hot streams' heat supply to the cold streams, i.e.\ an amount $q_{i,s,j,t}$ of heat exchanged between hot stream $i\in H$ in temperature interval $s\in T$ and cold stream $j\in C$ in temperature interval $t\in T$.
Heat may only flow to the same or a lower temperature interval, i.e.\ $q_{i,s,j,t}=0$, for each $i\in H$, $j\in C$ and $s,t\in T$ such that $s>t$.
A hot stream $i\in H$ and a cold stream $j\in C$ are \emph{matched}, if there is a positive amount of heat exchanged between them, i.e.\ $\sum_{s,t\in T}q_{i,s,j,t}>0$.
The objective is to find a feasible solution minimizing the number of matches $(i,j)$.

\begin{figure*}[t!]
\centering

\begin{subfigure}[t]{0.45\textwidth}
\centering
\includegraphics[]{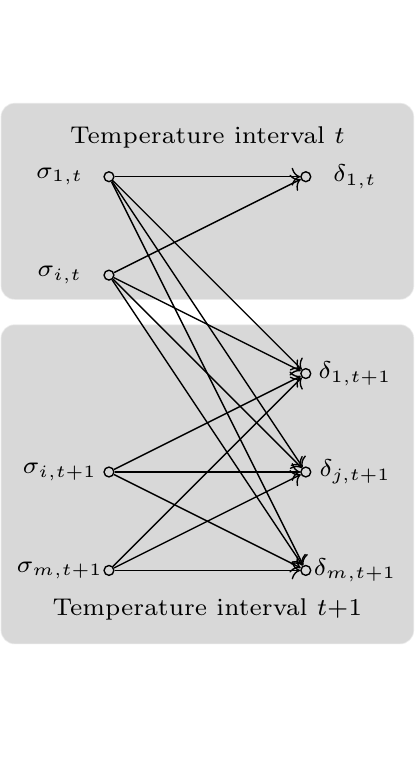}
\vspace*{-1cm}
\caption{ Transportation Model}
\label{Fig:transportation}
\end{subfigure}
\begin{subfigure}[t]{0.45\textwidth}
\centering
\includegraphics[]{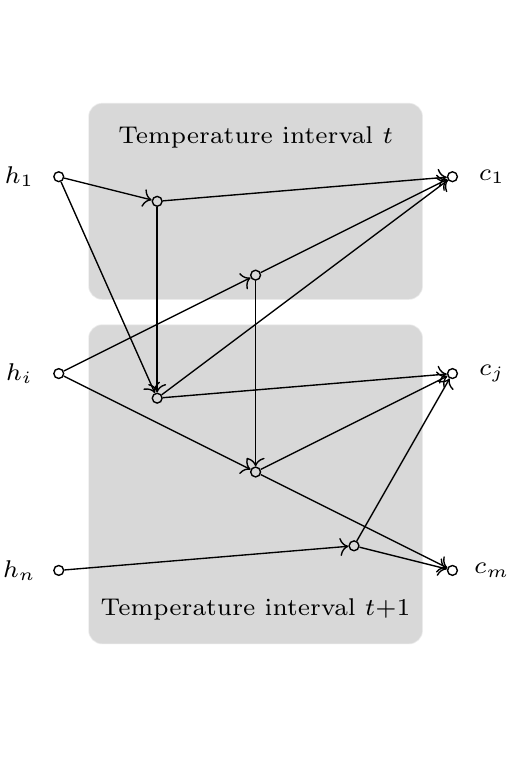}
\vspace*{-1cm}
\caption{Transshipment Model}
\label{Fig:transshipment}
\end{subfigure}
\vspace*{-0.5cm}
\caption{
In the transportation model \citep{Cerda1983matches}, each hot stream $i$ supplies $\sigma_{i,t}$ units of heat in temperature interval $t$ which can be received, in the same or a lower temperature interval, by a cold stream $j$ which demands $\delta_{j,t}$ units of heat in $t$. 
In the transshipment model \citep{Papoulias1983}, there are also intermediate nodes transferring residual heat to a lower temperature interval.
This figure is adapted from \citet{Furman2004}.
}
\label{Fig:Matches_Models}
\end{figure*}

\subsection{Mathematical Models}

This section presents a quadratic programming (QP) formulation for the multistage minimum utility cost problem and an MILP formulation for the minimum number of matches problem.

\subsubsection{Multistage Minimum Utility Cost Problem}

In the Eq.\ (\ref{Eq:Multistage}) QP formulation, continuous variables $Q_i^{CU}$ and $Q_j^{HU}$ compute the heat transferred from hot stream $i\in H$ to the cold utility and from the hot utility to cold stream $j\in C$.
Continuous variables $t_{i,k}$ and $t_{j,k}$ correspond to the temperature of hot stream $i\in H$ and cold stream $j\in C$ when exiting and entering stage $k\in S$, respectively.
Continuous variables $t_{i,j,k}^H$ and $t_{i,j,k}^C$ express the exiting temperature of hot stream $i\in H$ and cold stream $j\in C$ in heat exchanger $(i,j,k)$, respectively.
Continuous variables $f_{i,j,k}^H,f_{i,j,k}^C$ model the flow rate heat capacity of the hot and cold substream in heat exchanger $(i,j,k)$.
Auxiliary continuous variables $q_{i,j,k}$ are the heat exchanged via heat exchanger $(i,j,k)$.

{\allowdisplaybreaks {\footnotesize
\begin{subequations}
\label{Eq:Multistage}
\begin{align}
\min \quad & \sum_{i\in H}c^{CU}Q_i^{CU} + \sum_{j\in C}c^{HU}Q_j^{HU} \label{Eq:heat_cost} \\
& Q_i^{CU} = F_i(t_{i,\ell} - T_i^{\text{out}}) & i\in H \label{Eq:cold_utility} \\
& Q_j^{HU} = F_j(T_j^{\text{out}} - t_{j,0}) & j\in C \label{Eq:hot_utility} \\
& \sum_{j\in C} f_{i,j,k}^H = F_i^H & i\in H, k\in S \label{Eq:hot_split} \\
& \sum_{i\in H} f_{i,j,k}^C = F_j^C & j\in C, k\in S \label{Eq:cold_split} \\
& q_{i,j,k} = f_{i,j,k}^H (t_{i,k-1} - t_{i,j,k}^H) & i\in H, j\in C, k\in S \label{Eq:hot_heat} \\
& q_{i,j,k} = f_{i,j,k}^C (t_{i,j,k}^C - t_{j,k}) & i\in H, j\in C, k\in S \label{Eq:cold_heat} \\
& F_it_{i,k} = \sum_{j\in C} f_{i,j,k}^H t_{i,j,k}^H & i\in H, k\in S \label{Eq:hot_mixing} \\
& F_jt_{j,k-1} = \sum_{i\in H} f_{i,j,k}^C t_{i,j,k}^C & j\in C, k\in S \label{Eq:cold_mixing} \\
& t_{i,k-1} \leq t_{i,k} & i\in H, k\in S \label{Eq:hot_motonicity} \\
& t_{j,k-1} \leq t_{j,k} & j\in C, k\in S \label{Eq:cold_motonicity} \\
& T_i^{\text{in}} = t_{i,0} \geq t_{i,\ell} \geq T_i^{\text{out}} & i\in H \label{Eq:hot_bounds} \\
& T_j^{\text{out}} \geq t_{j,0} \geq t_{j,\ell} = T_j^{\text{in}} & j\in H \label{Eq:cold_bounds} \\
& t_{i,k}, t_{j,k}, t_{i,j,k}^H, t_{i,j,k}^C \geq 0 & i\in H, j\in C, k\in S \label{Eq:temperature_positive} \\ 
& f_{i,j,k}^H, f_{i,j,k}^C, q_{i,j,k}, Q_i^{CU}, Q_j^{HU}\geq 0 & i\in H, j\in C, k\in S \label{Eq:heat_positive}
\end{align}
\end{subequations}
}}
Expression (\ref{Eq:heat_cost}) minimizes the total heating utility cost.
Constraints (\ref{Eq:cold_utility}) and (\ref{Eq:hot_utility}) compute the heat absorbed by cold utilities and the heat supplied by hot utilities.
Constraints (\ref{Eq:hot_split}) and (\ref{Eq:cold_split}) divide the flow rate heat capacity of each stream fractionally to its corresponding substreams.
Constraints (\ref{Eq:hot_heat}) and (\ref{Eq:cold_heat}) compute the heat load exchanged between each pair of streams and enforce heat conservation.
Constraints (\ref{Eq:hot_mixing}) and (\ref{Eq:cold_mixing}) compute temperature of each stream by mixing substreams.
Constraints (\ref{Eq:hot_motonicity}) and (\ref{Eq:cold_motonicity}) enforce temperature monotonicity.
Constraints (\ref{Eq:hot_bounds}) and (\ref{Eq:cold_bounds}) assign initial temperature values and impose final temperature bounds.
Finally, Constraints (\ref{Eq:temperature_positive}) and (\ref{Eq:heat_positive}) ensure that all variables are non-negative.

\subsubsection{Minimum Number of Matches Problem}

The minimum number of matches can be formulated as an MILP using either the transportation, or the transshipment model in Figure \ref{Fig:Matches_Models}.
The former model represents heat as a commodity transported from supply nodes to destination nodes.
For each hot stream $i\in H$, there is a set of supply nodes, one for each temperature interval $s\in T$ with $\sigma_{i,s}>0$.
For each cold stream $j\in C$, there is a set of demand nodes, one for each temperature interval $t\in T$ with $\delta_{j,t}>0$.
There is an arc between the supply node $(i,s)$ and the destination node $(j,t)$ if $s\leq t$, for each $i\in H$, $j\in C$ and $s,t\in T$.
Continuous variable $q_{i,s,j,t}$ specifies the heat transferred from hot stream $i\in H$ in temperature interval $s\in T$ to cold stream $j\in C$ in temperature interval $t\in T$.
Binary variable $y_{i,j}$ indicates whether streams $i\in H$ and $j\in C$ are matched.
Big-M parameter $U_{i,j}$ bounds the amount of heat exchanged between every pair of hot stream $i\in H$ and cold stream $j\in C$, e.g.\ $U_{i,j}=\min\{h_i,c_j\}$.
Then, the problem can be modeled with formulation (\ref{Eq:Matches}).

{\footnotesize
{\allowdisplaybreaks
\begin{subequations}
\label{Eq:Matches}
\begin{align}
\text{min} & \sum_{i \in H}\sum_{j \in C} y_{i,j} \label{TransportationMIP_Eq:ObjMinMatches} \\ 
& \sum_{j\in C}\sum_{t\in T} q_{i,s,j,t} = \sigma_{i,s} & i\in H, s\in T \label{TransportationMIP_Eq:HotStreamConservation}\\
& \sum_{i\in H}\sum_{s\in T} q_{i,s,j,t} = \delta_{j,t} & j\in C, t\in T \label{TransportationMIP_Eq:ColdStreamConservation}\\
& \sum_{s,t\in T} q_{i,s,j,t}\leq U_{i,j}\cdot y_{i,j} & i\in H, j\in C \label{TransportationMIP_Eq:BigM_Constraint}\\
& q_{i,s,j,t} = 0 & i\in H, j\in C, s,t\in T: s> t \label{TransportationMIP_Eq:ThermoConstraint} \\
& y_{i,j} \in \{0,\,1\},q_{i,s,j,t}\geq 0 & i\in H, j\in C,\; s,t\in T \label{TransportationMIP_Eq:IntegralityConstraints}
\end{align}
\end{subequations}
}
}
Expression (\ref{TransportationMIP_Eq:ObjMinMatches}), the objective function, minimizes the number of matches.
Equations (\ref{TransportationMIP_Eq:HotStreamConservation}) and (\ref{TransportationMIP_Eq:ColdStreamConservation}) ensure heat conservation.
Equations (\ref{TransportationMIP_Eq:BigM_Constraint}) enforce a match between a hot and a cold stream if they exchange a positive amount of heat.
Equations (\ref{TransportationMIP_Eq:BigM_Constraint}) are \emph{big-M constraints}.
Equations (\ref{TransportationMIP_Eq:ThermoConstraint}) ensure that no heat flows to a hotter temperature.

\subsection{Computational Complexity and Approximation Algorithms}

\citet{Furman2001} show that minimum number of matches problem is strongly $\mathcal{NP}$-hard, even in the special case with a single temperature interval, through a reduction from 3-Partition \citep{Garey2002}.
\citet{Letsios2018} present an $\mathcal{NP}$-hardness reduction from bin packing.
\citet{Furman2001} demonstrate that the more general hyperstructure, multistage, and sequential heat exchanger network synthesis are all strongly $\mathcal{NP}$-hard as they can be reduced to the minimum number of matches problem.
On the positive side, the minimum utility cost problem in sequential heat exchanger network synthesis can be formulated as an LP and is, therefore, polynomially solvable.
The complexity of the multistage minimum utility cost problem is an intriguing open question.


\citet{Furman2004} initiate the design of approximation algorithms for heat exchanger network synthesis problems.
In particular, they investigate the approximability of the minimum number of matches problem and propose (i) a collection of greedy and relaxation rounding heuristics, (ii) an $O(r)$-approximation algorithm, where $r$ is the number of temperature intervals, and (iii) a 2-approximation ratio for the single temperature interval subproblem.
\citet{Letsios2018} classify the heuristics for the minimum number of matches of problem into relaxation rounding, water filling, and greedy packing.
For the general problem, they show (i) an $\Omega(n)$ bound on the approximation ratio of deterministic LP rounding, (ii) an $\Omega(k)$ bound on the approximation ratio of greedy water filling, and (iii) a positive $O(\log n/\epsilon)$ ratio for greedy packing.
For the single temperature interval subproblem, they propose an improved $1.5$-approximation algorithm.

\section{Concluding Remarks and Future Directions}
\label{Section:Conclusion}

This paper discusses ways of using approximation algorithms for solving challenging PSE problems and reports state-of-the-art examples motivating this line of work.
We outline applications in: (i) mathematical modeling, (ii) problem classification, (iii) design of solution methods, and (iv) dealing with uncertainty.
In order to exemplify the proposed investigations, we consider three fundamental PSE optimization problems: pooling, process scheduling, and heat exchanger network synthesis. There are many other possible PSE applications, e.g.\ in at the intersection between scheduling and control \citep{pistikopoulos2016towards,DIAS2018139,DAOUTIDIS2018179,DIAS2019, etesami2019unifying,TSAY201922}, which provide additional and interesting challenges.

This paper presents formal problem descriptions, standard mathematical programming formulations, brief literature surveys, and prepares the ground for investigating three fundamental PSE optimization  problems from an approximation algorithms perspective.
Some future challenges we see in this area are as follows:


\begin{enumerate}
\item Pooling remains $\mathcal{NP}$-hard when each raw material supply, final product demand, and quality attribute must be equal to a fixed value.
In these fixed-value cases, pooling is a variant of standard multicommodity flow problems, which are among the most extensively studied combinatorial objects in TCS.
Extensions of the well-known min-cut max-flow theorem to the multicommodity flow setting result in tight relaxations and dual multicut bounds 
\citep{Garg1996,Leighton1999}. \\
\emph{
Can we derive strong algorithms for large-scale instances via connections to multicommodity flow?}


\item Pooling becomes more tractable in the case of sparse instances. 
Furthermore, discretization enables efficient pooling solving with exact methods. \\
\emph{Using the quality of sparse and discrete relaxations, can we compute problem classifications to develop useful trade-offs between solution quality and running time efficiency?}

\item Process scheduling involves tasks with variable processing times to determine the batch sizes.
Scheduling with controllable processing times is an active operations research area dealing with this setting \citep{Shabtay2007,Shioura2018}.
In TCS, analogous investigations have taken place in the context of speed scaling where a processing unit may modify its speed to save energy and task processing times are decision variables
\citep{Albers2010,Albers2017,Angel2019,Bampis2015,Bampis2016,Bampis2018,Bansal2007,Yao1995}. \\
\emph{Can we apply techniques for obtaining algorithms with analytically proven performance guarantees, including network flows, convex relaxations, and submodular optimization for solving PSE problem instances?}

\item State-task network problems are strongly related to precedence-constrained, shop, and resource-constrained project scheduling
\citep{Bampis2014, Hall1989, Kone2011}. \\
\emph{Could the different relaxations developed for these scheduling variants result in stronger mathematical modeling strategies for process scheduling?}

\item Determining the computational complexity of the multistage minimum utility cost problem is an intriguing future direction.
Because of stream mixing, the problem exhibits commonalities with pooling.
However, no hardness reduction formalizes this insight of domain experts. \\
\emph{Could efficient approximation algorithms for the multistage minimum utility cost problem assist in solving simultaneous heat exchanger network synthesis at industrial scales?}

\item The minimum number of matches problem remains a major bottleneck in heat exchanger network synthesis.
The problem can be considered as a special two-dimensional packing where the vertical and horizontal axis correspond to temperature and flow rate heat capacity, respectively. \\ 
\emph{Could we take advantage of this packing nature to derive stronger formulations?}
\end{enumerate}


\section*{Acknowledgements}

\noindent
This work was funded by Engineering \& Physical Sciences Research Council (EPSRC) grant EP/M028240/1, an EPSRC DTP (award ref. 1675949) to RBL, an EPSRC Center for Doctoral Training in High Performance Embedded and Distributed Systems (EP/L016796/1) studentship to MM, an EPSRC/Schlumberger CASE studentship to JW (EP/R511961/1, voucher 17000145), and an EPSRC Research Fellowship to RM (grant number EP/P016871/1).

\singlespacing

\bibliographystyle{elsarticle-harv}\biboptions{authoryear}
\bibliography{refs_CACE}


\appendix

\section{Nomenclature}

\subsection{Pooling Problem}\label{app:pooling_notation}

\begin{longtable}{l l l}
        \toprule
        Type & Name & Description\\
         \midrule
         Sets & $I$ & Inputs nodes, raw materials \\
         & $L$ & Pool nodes, intermediate products \\
         & $J$ & Output nodes, end products \\
         & $X$ & Input-to-pool arcs\\
         & $Y$ & Pool-to-output arcs\\
         & $Z$ & Input-to-output arcs\\
         & $K$ & Quality attributes \\
        \midrule 
        Indices &  $i$ & Input \\
         & $l$ & Pool \\
         & $j$ & Output \\
         & $k$ & Attribute \\
         \midrule
         Parameters & $\ci$ & Raw material unitary cost \\
         & $d_j$ & End product unitary profit \\
         & $\Ail,\Aiu$ & Raw material supply bounds \\
         & $\Sl$ & Pool capacity\\
         & $\Djl,\Dju$ & End product demand bounds \\
         & $\Cik$ & Raw material quality attribute \\
         & $\Pjkl,\Pjku$ & End product quality attribute range \\
         & $\Delta_i^{\stout},\Delta_l^{\stout}$ & Input, pool node out-degree \\
         & $\Delta_l^{\stin},\Delta_j^{\stin}$ & Pool, output node in-degree \\
         \midrule
         Variables & $\xil$ & Input-to-pool flow\\
         & $\ylj$ & Pool-to-output flow\\
         & $\zij$ & Bypass input-to-output flow \\
         & $\vilj$ & Path flow \\
         & $\qil$ & Input-to-pool fractional flow \\ 
         & $\plk$ & Intermediate product quality attribute \\       
\bottomrule
\label{tab:pp-formulations}
\end{longtable}

\subsection{Process Scheduling}\label{app:scheduling_notation}

\begin{longtable}{l l l}
        \toprule
        Type & Name & Description\\
        \midrule
         Sets & $S$ & States \\
              & $I$ & Tasks \\
              & $A$ & State-task network arcs \\
              & $A^+,A^-$ & Production, consumption arcs \\
              & $J$ & Units\\
              & $T$ & Time slots, time points\\
              & $I_j$ & Tasks that unit  $j$ may execute \\ 
              & $J_i$ & Units capable of performing task $i$ \\
              & $S_i^-, S_i^+$ & States consumed, produced by task $i$ \\
              & $I_s^-, I_s^+$ & Tasks consuming, producing state $s$\\      
        \midrule
         Indices    & $s$ & State\\
                    & $i$ & Task\\
                    & $j$ & Unit\\
                    & $t,t'$ & Time, time slot\\
                    & $k$ & Time point \\
         \midrule
         Parameters & $\pij$ & Processing time of task $i$ on unit $j$\\
                    & $\bijl,\biju$ & Minimum, maximum processing capacity of unit $j$ for task $i$\\
                    & $f_{i,s}^-,f_{i,s}^+$ & Material fraction entering, exiting as state $s$ for task $i$\\
                    & $\qism$ & Fraction of material for task $i$ entering as state $s$\\                    
                    & $d_s$  & Demand for state $s$ \\
                    & $\ell$ & Time horizon, number of time slots or time points \\
        \midrule
         Variables & $z$ & Makespan \\ 
                   & $\xijt$ & Indicates whether task $i$ begins processing on unit $j$ at time $t$ \\
                   & $\bijt$ & Batch size of task $i$ starting on unit $j$ at time $t$ \\     
                   & $\yst$  & Stored amount of state $s$ at time $t$ \\    
                   & $t_k$ & Time point $k$ \\   
                   & $x_{i,k}^S, x_{i,k}^F$ & Indicate whether task $i$ begins, finishes at time point $k$ \\
                   & $t_{i,k}^S, t_{i,k}^F$ & Starting, finishing time of task $i$ starting at time point $k$ \\  
                   & $p_{i,k}$ & Processing time of task $i$ starting at time point $k$ \\
                   & $b_{i,k}$ & Batch size of task $i$ starting at time point $k$ \\     
                   & $y_{s,k}$ & Amount of state $s$ at time point $k$ \\
         \bottomrule
\label{tab:notation-stn}
\end{longtable}

\subsection{Heat Exchanger Network Synthesis}\label{app:hens_notation}

\begin{longtable}{l l l}
    \toprule
        Type & Name & Description \\
        \midrule
        Sets & $H$ & Hot streams \\
        & $C$ & Cold streams \\
        & $S$ & Stages \\
        & $T$ & Temperature intervals \\
        \midrule
        Indices & $i$ & Hot stream \\
        & $j$ & Cold stream \\
        & $k$ & Stage \\
        & $s$,$t$ & Temperature interval \\
        \midrule
        Parameters & $T^{\stin}_i, T^{\stout}_i$ & Hot stream $i$ inlet and outlet temperature \\
        & $T^{\stin}_j, T^{\stout}_j$ & Cold stream $j$ inlet and outlet temperature \\ 
        & $F_i, F_j$ & Flow rate heat capacity of hot stream $i$ and cold stream $j$ \\
        & $c_{HU}, c_{CU}$ & Heating and cooling utility unitary cost \\
        & $\ell$ & Number of stages \\
        & $r$ & Number of temperature intervals \\
        & $T_t$ & $t$-th greatest discrete inlet / outlet temperature value \\
        & $h_i$ & Heat load of hot stream $i$ \\
        & $c_j$ & Heat load of cold stream $j$ \\
        & $\sigma_{i,s}$ & Heat supply of hot stream $i$ at temperature interval $s$ \\
        & $\delta_{j,t}$ & Heat demand of cold stream $j$ at temperature interval $t$ \\
        & $U_{i,j}$ & Upper bound on heat exchanged between hot stream $i$ and cold stream $j$ \\
        & $\Delta T_{\min}$ & Minimum heat approach temperature \\
        \midrule
        Variables & $t_{i,k}$ & Final temperature of hot stream $i$ at stage $k$ \\
        & $t_{j,k}$ & Initial temperature of cold stream $j$ at stage $k$ \\
        & $t_{i,j,k}^H, t_{i,j,k}^C$ & Temperature of heat exchanger $(i,j,k)$ in the hot and cold side \\
        & $f_{i,j,k}^H, f_{i,j,k}^C$ & Flow rate heat capacity of heat exchanger $(i,j,k)$ in the hot and cold side \\
        & $q_{i,j,k}$ & Heat transferred via heat exchanger $(i,j,k)$ \\
        & $Q_i^{CU}$ & Cold utility heat load from hot stream $i$ \\
        & $Q_j^{HU}$ & Hot utility heat load to cold stream $i$ \\
        & $y_{i,j,}$ & Binary indicating whether hot stream $i$ is matched with cold stream $j$ \\
        & $q_{i,s,j,t}$ & Heat exchanged between hot stream $i$ at temperature interval $s$ and cold \\
        & & stream $j$ at temperature interval $t$ \\
        \bottomrule
\label{tab:hens-notation}
\end{longtable}

\end{document}